\documentclass{article}

\usepackage{arxiv}

\usepackage[utf8]{inputenc} % allow utf-8 input
\usepackage[T1]{fontenc}    % use 8-bit T1 fonts
\usepackage{hyperref}       % hyperlinks
\usepackage{url}            % simple URL typesetting
\usepackage{booktabs}       % professional-quality tables
\usepackage{amsfonts}       % blackboard math symbols
\usepackage{nicefrac}       % compact symbols for 1/2, etc.
\usepackage{microtype}      % microtypography
\usepackage{lipsum}		% Can be removed after putting your text content
\usepackage{graphicx}
\usepackage{natbib}
\usepackage{doi}
\usepackage{amsmath, amsthm, amssymb}
\newcommand{\seq}[1]{\left<#1\right>}
\newcommand{\norm}[1]{\left\Vert#1\right\Vert}
\newcommand{\abs}[1]{\left\vert#1\right\vert}
\newcommand{\set}[1]{\left\{#1\right\}}
\newcommand{\bra}[1]{\left(#1\right)}
\newcommand{\bras}[1]{\left(#1\right]}
\newcommand{\brass}[1]{\left[#1\right)}
\newcommand{\sbra}[1]{\left[#1\right]}
\newcommand{\vertiii}[1]{{\left\vert\kern-0.25ex\left\vert\kern-0.25ex\left\vert #1
    \right\vert\kern-0.25ex\right\vert\kern-0.25ex\right\vert}}
\newcommand{\R}{\mathbb R}
\newcommand{\N}{\mathbb N}
\renewcommand{\to}{\longrightarrow}
\newcommand{\h}{\mathcal{H}}
\newcommand{\bh}{\mathcal{B}(\mathcal{H})}

\renewcommand{\c}{\mathbb{C}}
\newcommand{\diag}{\mathrm{diag}}

\newtheorem{theorem}{Theorem}[section]
\newtheorem{lemma}[theorem]{Lemma}

\newtheorem{corollary}[theorem]{Corollary}
\newtheorem{remark}[theorem]{Remark}
\newcommand\mystyle{\everymath{\displaystyle}}
\mystyle
\title{New Improvements to Heron and Heinz Inequality Using Matrix Techniques}

\author{\href{https://orcid.org/0000-0002-3816-5287}{\includegraphics[scale=0.06]{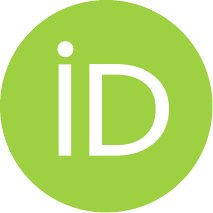}\hspace{1mm}M.H.M.~Rashid}\thanks{Corresponding Author} \\
	Department of Mathematics\&Statistics\\Faculty of Science P.O.Box(7)\\
	Mu'tah University University\\
	Mu'tah-Jordan \\
	\texttt{mrash@mutah.edu.jo}
	 \AND
	{Wael Mahmoud Mohammad Salameh}\\
	{Faculty of Information Technology}\\
 {Abu Dhabi University, Abu Dhabi 59911}\\{United Arab Emirates}\\ 
	{wael.salameh1@adu.ac.ae}\\
}

\renewcommand{\shorttitle}{\textit{arXiv} Template}

\hypersetup{
pdftitle={New Improvements to Heron and Heinz Inequality Using Matrix Techniques},
pdfsubject={q-bio.NC, q-bio.QM},
pdfauthor={M.H.M.Rashid, Wael Mahmoud Mohammad Salameh},
pdfkeywords={Heinz mean inequalities, positive semi-definite matrices, Hilbert-Schmidt norm, Young inequality},
}

\begin{document}
\maketitle

\begin{abstract}
	This paper undertakes a thorough investigation of matrix means interpolation and comparison. We expand the parameter $\vartheta$ beyond the closed interval $[0,1]$ to cover the entire positive real line, denoted as $\mathbb{R}^+$. Furthermore, we explore additional outcomes related to Heinz means. We introduce new scalar adaptations of Heinz inequalities, incorporating Kantorovich's constant, and enhance the operator version. Finally, we unveil refined Young's type inequalities designed specifically for traces, determinants, and norms of positive semi-definite matrices.
\end{abstract}

% keywords can be removed

\section{Introduction}
Consider the algebra of complex matrices of size $n\times n$, denoted as $M_n(\mathbb{C})$. A matrix $T$ in $M_n(\mathbb{C})$ is considered positive semi-definite, written as $T\geq 0$, if it is Hermitian and satisfies $\langle Tx,x\rangle\geq 0$ for all vectors $x$ in $\mathbb{C}^n$. If, for a Hermitian matrix $T$ in $M_n(\mathbb{C})$, $\langle Tx,x\rangle>0$ holds for all nonzero vectors $x$ in $\mathbb{C}^n$, it is termed a positive definite matrix, denoted as $T>0$. The set of all positive matrices is denoted as $M_n^+(\mathbb{C})$, and the subset of definite matrices within $M_n^+(\mathbb{C})$ is represented as $M_n^{++}(\mathbb{C})$. The Schur product of two matrices $T=[t_{ij}]_{i,j}$ and $S=[s_{ij}]_{i,j}$ in $M_n(\mathbb{C})$ is defined as the matrix $T\circ S$ with entries $t_{ij}s_{ij}$.
 A norm $\vertiii{.}$ on the set of complex matrices of size $n\times n$, denoted as $M_n(\mathbb{C})$, is termed unitarily invariant if $\vertiii{UAV}=\vertiii{T}$ for any matrix $T$ in $M_n(\mathbb{C})$ and for all unitary matrices $U$ and $V$ in $M_n(\mathbb{C})$.

For a matrix $T=[t_{ij}]\in M_n(\mathbb{C})$, the Hilbert-Schmidt norm (also known as the Frobenius norm) and the trace norm of $T$ are defined as follows
\begin{equation}\label{Sch1}
  \norm{T}_2=\bra{\sum_{j=1}^{n}s_j^{2}(T)}^{\frac{1}{2}},\qquad \norm{T}_1=tr(|T|)=\sum_{j=1}^{n}s_j(T)
\end{equation}
Here, $s_1(T) \geq s_2(T)\geq \cdots\geq s_n(T)\geq 0$ represent the singular values of $T$, which are the eigenvalues of the positive matrix $|T|=\sqrt{T^{*}T}$ arranged in decreasing order and repeated according to multiplicity. The symbol $tr(.)$ denotes the usual trace operation.

It's important to note that the mathematical norms $\norm{\cdot}_2$ and $\norm{\cdot}_1$ are widely recognized for being unitarily invariant.

The classic Young's inequality for non-negative real numbers states that if  $\rho,\sigma\geq 0$ and $0\leq \kappa\leq 1$, then
\begin{equation}\label{Y1}
  \rho^{\kappa}\sigma^{1-\kappa}\leq \kappa \rho+(1-\kappa)\sigma
\end{equation}
Equality occurs if and only if $\rho=\sigma$. When $\kappa$ is $\frac{1}{2}$, substituting into the inequality yields the arithmetic-geometric mean inequality
\begin{equation}\label{Y2}
  \sqrt{\rho\sigma}\leq \frac{\rho+\sigma}{2}.
\end{equation}
Manasarah and Kittaneh, as presented in \cite{KM}, improved Young's inequality with the following refinement
\begin{equation}\label{Y3}
  \bra{\rho^{\kappa}\sigma^{1-\kappa}}^{m}+r_0^m\bra{\rho^{\frac{m}{2}}-\sigma^{\frac{m}{2}}}^2\leq \bra{\kappa \rho^{r}+(1-\kappa)\sigma^{r}}^{\frac{m}{r}},\quad r\geq1
\end{equation}
where $m\in\N$ and $r_0=\min\{\kappa,1-\kappa\}$.\\
%%%%%%%%%%%%%%%%%%%%%%%%%%%%%%%%%%%%%%%%%%%%%%%%%%%%%%%%%%%%%%%%%%%%%%%%%%%%%%%%%%%%
%========================================================================
\indent The Kantorovich constant, denoted as $K(t,2)$, is defined as $\frac{(t+1)^2}{4t}$. It possesses several key properties: $K(1,2)=1$, $K(t,2)=K\bra{\frac{1}{t},2}\geq 1$ $(t>0)$ and $K(t, 2)$ is
monotone increasing on $[1,\infty)$, and monotone decreasing on $(0,1]$. For more detailed information about the Kantorovich constant, interested readers can refer to \cite{NL, Rash-Feras, ZSF, ZW1}..\\
%=====================================================================
The following multiplicative refinement and reversal of Young's inequality, expressed in terms of Kantorovich's constant, can be stated as follows
\begin{equation}\label{A5}
  K(h,2)^r\rho\sharp_{\kappa}\sigma\leq \rho\nabla_{\kappa}\sigma\leq K(h,2)^{R}\rho\sharp_{\kappa}\sigma,
\end{equation}
where $\rho$ and $\sigma$ are both greater than 0, $\kappa$ belongs to the interval $[0,1]$, $r$ is the minimum of $\kappa$ and $1-\kappa$, $R$ is the maximum of $\kappa$ and $1-\kappa$, and $h$ is defined as $\frac{\sigma}{\rho}$.

The second inequality in (\ref{A5}) is credited to Liao et al. \cite{LWZ}, while the first one is attributed to Zou et al. \cite{ZSF}.
%============================================================
In \cite{WZ}, the authors obtained another improvement of the Young inequality
and its reverse as follows:
\begin{equation}\label{A6}
 r(\sqrt{\rho}-\sqrt{\sigma})^2+ K(\sqrt{h},2)^{r'}\rho\sharp_{\kappa}\sigma\leq \rho\nabla_{\kappa}\sigma,
\end{equation}
and
\begin{equation}\label{A7}
  \rho\nabla_{\kappa}\sigma\leq K(\sqrt{h},2)^{-r'}\rho\sharp_{\kappa}\sigma+R(\sqrt{\rho}-\sqrt{\sigma})^2
\end{equation}
where $h=\frac{\sigma}{\rho}$, $r=\min\{\kappa,1-\kappa\}$, $R=\max\{\kappa,1-\kappa\}$ and $r'=\min\{2r,1-2r\}$.
In addition, another kind of the reversal of Young inequality utilizing Kantorovich's constant is described
 in \cite{LWZ} with the same notation as above.
\begin{equation}\label{A8}
  \rho\nabla_{\kappa}\sigma-R(\sqrt{\rho}-\sqrt{\sigma})^2\leq K(\sqrt{h},2)^{R'}\rho\sharp_{\kappa}\sigma,
\end{equation}
where $R'=\max\{2r,1-2r\}$.\\
%=============================================================================
For $\kappa$ in the range of $[0,1]$ and two non-negative real numbers $\rho$ and $\sigma$, the Heinz mean serves as an interpolation between the $\kappa$-arithmetic mean and the $\kappa$-geometric mean. These are defined by the expression
\begin{equation}\label{A18}
H_{\kappa}(\rho,\sigma)=\frac{\rho\sharp_{\kappa} \sigma+\rho\sharp_{1-\kappa}\sigma}{2},
\end{equation}
where $\rho\sharp_{\kappa} \sigma=\rho^{\kappa}\sigma^{1-\kappa}$ represents the $\kappa$-geometric mean.

The Heinz mean possesses certain properties, including convexity concerning $\kappa$ within the interval $[0,1]$. Its minimum occurs at $\kappa=\frac{1}{2}$, and its maximum values are found at $\kappa=0$ and $\kappa=1$. Additionally, the following inequalities are true
\begin{equation}
\sqrt{\rho\sigma}\leq H_{\kappa}(\rho,\sigma)\leq \frac{\rho+\sigma}{2}.
\end{equation}\label{A9}
It is worth noting that the function $H_{\kappa}(\rho,\sigma)$ exhibits symmetry with respect to the point $\kappa=\frac{1}{2}$, meaning that $H_{\kappa}(\rho,\sigma)=H_{1-\kappa}(\rho,\sigma)$.\\
The Heron mean is defined by the expression
\begin{equation}\label{A15}
  F_{\vartheta}(\rho,\sigma)=(1-\vartheta)\sqrt{\rho\sigma}+\vartheta\bra{\frac{\rho+\sigma}{2}},\,\,\vartheta\in [0,1]\,\,\mbox{and}\,\, \rho,\sigma\in\R^+.
\end{equation}
where $\vartheta$ takes values in the interval $[0,1]$, and $\rho$ and $\sigma$ are positive real numbers.

Evidently, the Heron mean serves as a linear interpolation between the arithmetic and geometric means. It adheres to the inequality $F_{\vartheta}\leq F_{\varrho}$ whenever $\vartheta\leq \varrho$, with both $\vartheta$ and $\varrho$ belonging to the positive real numbers.

In a study by Bhatia published in \cite{Bhatia}, it was demonstrated that for $\vartheta(\kappa)=(2\kappa-1)^2$ and $\kappa$ within the range of $[0,1]$, the following relation holds
\begin{equation}\label{A20}
  H_{\kappa}(\rho,\sigma)\leq F_{\vartheta(\kappa)}(\rho,\sigma).
\end{equation}
%======================================================================================
\indent Our paper is structured as follows: In the upcoming section, we will conduct an in-depth investigation into matrix interpolation and mean comparisons. This analysis extends the scope of $\vartheta$ beyond the closed interval $[0,1]$ to include all positive real numbers, represented as $\mathbb{R}^+$. Additionally, we will explore additional findings pertaining to Heinz means. Section three is dedicated to exploring refinements in Heinz inequality, incorporating the Kantorovich constant. Section 4 focuses on examining enhanced variations of Heinz-type operator inequalities and their corresponding reversals. Finally, in section 5, we present refined inequalities of Young's type, specifically designed for traces, determinants, and norms of positive semi-definite matrices.
%===============================================================================
\section{Full Interpolation of Matrix Variants of Heron and Heinz  MEANS}
%=========================================================================
In the paper referenced as \cite{Bhatia}, R. Bhatia established a noteworthy result. In particular, it was shown that for values of $\vartheta$ in the interval $[0,1/2]$, the function $\psi(\vartheta)$ adheres to the inequality $\psi(\vartheta) \leq \psi(1/2)$. Here, $\psi(\vartheta)$ denotes one of the potential matrix formulations of Equation (\ref{A15}), and its definition is as follows
\begin{equation}\label{BA1}
\psi(\vartheta)=\vertiii{(1-\vartheta)T^{1/2}XS^{1/2}+\vartheta\bra{\frac{TX+XS}{2}}},
\end{equation}
This definition involves matrices $T$, $S$, and $X$, subject to the conditions that $T$ and $S$ belong to the set of positive definite matrices in $\mathbb{C}$, denoted as $M_n^{++}(\mathbb{C})$, and $X$ is a member of the set of $n\times n$ matrices over $\mathbb{C}$, denoted as $M_n(\mathbb{C})$.

For further insights into the matrix formulations of Equation (\ref{A18}) and Equation (\ref{A20}), as well as additional details, interested readers are encouraged to refer to the following references: \cite{Bhatia}, \cite{Horn}, \cite{HJ}, \cite{HK}, and \cite{KS}.

Within the context of this article, the author endeavors to demonstrate that for $\vartheta$ values within the interval $[0,1/2]$, the function $\psi(\vartheta,\kappa)$ adheres to the inequality $\psi(\vartheta,\kappa) \leq \psi(1/2,\kappa)$. Additionally, it is asserted that $\psi(\vartheta,\kappa)$ displays an increasing trend as $\vartheta$ varies within the range $[1/2,\infty)$.

This result serves as a generalization of the previously established monotonic property associated with the matrix version of Equation (\ref{A15}). This generalization mirrors the behavior of $F_{\vartheta}(\rho,\sigma)$ for positive real numbers $a$ and $b$ when $\vartheta$ belongs to the set of positive real numbers, $\mathbb{R}^+$.

As a consequential outcome of these findings, the author will introduce a potential generalized matrix equivalent of Equation (\ref{A20}), which can be expressed as follows
\begin{equation}\label{BA2}
\frac{1}{2}\vertiii{T^{\mu}XS^{1-\mu}+T^{1-\mu}XS^{\mu}}\leq \vertiii{(1-\vartheta)T^{\kappa}XS^{1-\kappa}+\vartheta\bra{\frac{TX+XS}{2}}}
\end{equation}
This inequality is valid for particular values of $\mu\in [1/4,3/4]$, $\kappa\in [0,1]$, and $\vartheta\in [1/2,\infty)$.\\
%%%%%%%%%%%%%%%%%%%%%%%%%%%%%%%%%%%%%%%%%%%%%%%%%%%%%%%%%%%%%%%%%%%%%%%%%%%%%%%%%%%%%%%%%%%%%%%%%%%%%%%%
\indent In this section, we will undertake a thorough investigation of matrix interpolation and mean comparisons. This scrutiny broadens the range of $\vartheta$ from the closed interval $[0,1]$ to encompass the entirety of positive real numbers, denoted as $\mathbb{R}^+$. Furthermore, we will delve into additional findings associated with Heinz means.
%===================================Theorem2.1=======================================================
\begin{theorem}\cite{HJ}\label{HHM1} Let $T,S\in M_n(\c)$  such that $T$ is a positive semi-definite. Then
 $$\vertiii{T\circ S}\leq \max_{1\leq i\leq n} t_{ii}\vertiii{S},$$
where $t_{ii}$ for $i=1,2,\cdots,n$ are the diagonal entries of matrix $T$.
\end{theorem}
%===============================Lemma2.2========================================
\begin{lemma}\cite{Zhan}\label{HHM2}
  Let $\kappa_1,\kappa_2,\cdots,\kappa_n$ be positive numbers, $r\in [-1,1]$, and $t\in (-2,2]$. Then the $n\times n$ matrix
matrix $$\Gamma=\bra{\frac{\kappa_{i}^{r}+\kappa_{j}^{r}}{\kappa_{i}^{2}+t\kappa_i\kappa_j+\kappa_j^2}}$$
is positive semi-definite.
\end{lemma}
%==================================Theorem2.3=============================
\begin{theorem}\label{Rahma1} Let $T,S,X\in M_n(\c)$ such that $T$ and $S$ are positive semi-definite, $\kappa\in [0,1]$ and
$\vertiii{.}$ any unitarily invariant
norm, the function
$$\psi(\vartheta,\kappa)=\vertiii{(1-\vartheta)T^{\kappa}XS^{1-\kappa}+\vartheta\bra{\frac{TX+XS}{2}}}$$
is increasing for $\frac{1}{2}\leq\vartheta<\infty$ and $\psi(\vartheta,\kappa)\leq \psi\bra{\frac{1}{2},\kappa}$
for all $\vartheta\in \sbra{0,\frac{1}{2}}$.
\end{theorem}
\begin{proof} We first prove the result for $\vartheta>0$ and $T=S$, that is,
$$\vertiii{(1-\vartheta)T^{\kappa}XS^{1-\kappa}+\vartheta\bra{\frac{TX+XS}{2}}}=\frac{\vartheta}{2}Z(\vartheta),$$
where $Z(\vartheta)=\vertiii{Q(\vartheta)T^{\kappa}XT^{1-\kappa}+TX+XT}$ and $Q(\vartheta)=2\bra{\frac{1}{\vartheta}-1}$. We may assume without loss of generality,
$T=\diag\bra{\eta_1,\cdots,\eta_n}$, $\eta_j>0$. Then
\begin{eqnarray*}
  && Q(\vartheta)T^{\kappa}XT^{1-\kappa}+TX+XT=\bra{\bra{Q(\vartheta)\eta_i^{\kappa}\eta_j^{1-\kappa}+\eta_i+\eta_j}x_{ij}}_{i,j} \\
   &&= \bra{\frac{Q(\vartheta)\eta_i^{\kappa}\eta_j^{1-\kappa}+\eta_i+\eta_j}{Q(\varrho)\eta_i^{\kappa}\eta_j^{1-\kappa}+\eta_i+\eta_j}}_{i,j}
\circ\bra{Q(\varrho)T^{\kappa}XT^{1-\kappa}+TX+XT}\\
&&=E\circ \bra{Q(\varrho)T^{\kappa}XT^{1-\kappa}+TX+XT},
\end{eqnarray*}
where $E=\bra{\frac{Q(\vartheta)\eta_i^{\kappa}\eta_j^{1-\kappa}+\eta_i+\eta_j}{Q(\varrho)\eta_i^{\kappa}\eta_j^{1-\kappa}+\eta_i+\eta_j}}_{i,j}$. Now the matrix $E$ can be written as
$$\bra{1+\frac{(Q(\vartheta)-Q(\varrho))\eta_i^{\kappa}\eta_j^{1-\kappa}}{Q(\varrho)\eta_i^{\kappa}\eta_j^{1-\kappa}+\eta_i+\eta_j}}
=\bra{1}_{i,j}+\bra{\eta_i^{\kappa}\bra{\frac{Q(\vartheta)-Q(\varrho)}{Q(\varrho)\eta_i^{\kappa}\eta_j^{1-\kappa}+\eta_i+\eta_j}}\eta_j^{1-\kappa}}$$
which will be positive semidefinite if the matrix,
$$G=\bra{\frac{Q(\vartheta)-Q(\varrho)}{Q(\varrho)\eta_i^{\kappa}\eta_j^{1-\kappa}+\eta_i+\eta_j}}_{i,j}$$
is positive semidefinite. According to Lemma \ref{HHM2}, the latter matrix is positive semidefinite if and only if $Q(\vartheta)\geq Q(\varrho)$ and $Q(\varrho)\in \sbra{-2,2}$. Since $Q(\vartheta)=2\bra{\frac{1}{\vartheta}-1}$ is a continuous and decreasing function on the positive half-line, ranging from $[\frac{1}{2},\infty)$ into $[-2,2]$, it follows that $Q(\vartheta)\geq Q(\varrho)$ for all $\varrho\geq \vartheta$. Consequently, using Theorem \ref{HHM1}, we can deduce that $Z(\vartheta)\leq \bra{\frac{Q(\vartheta)+2}{Q(\varrho)+2}}T(\varrho)$. Thus, the result holds for $T=S$ and $\vartheta\geq \frac{1}{2}$.

For $\vartheta \in (0,1/2]$, we have $2\leq Q(\vartheta)<\infty$, and $Q(\vartheta)> Q\bra{\frac{1}{2}}=2$. Therefore, the matrix $E$ with $\varrho=\frac{1}{2} $ is positive semidefinite, as per Lemma \ref{HHM1}. The case $\vartheta=0$ is straightforward since, by Lemma \ref{HHM1}, the matrix
$$\bra{\frac{\eta_i^{\kappa}\eta_j^{1-\kappa}}{\eta_i^{\kappa}\eta_j^{1-\kappa}+\eta_i+\eta_j}}_{i,j}
=\bra{\eta_i^{\kappa}\bra{\frac{1}{\eta_i^{\kappa}\eta_j^{1-\kappa}+\eta_i+\eta_j}}\eta_j^{1-\kappa}}_{i,j}$$
is positive semidefinite. Thus, we have established the desired result for this case, i.e., $\vartheta Z(\vartheta)\leq \frac{1}{2}Z\bra{\frac{1}{2}}$. In other words, $\psi(\vartheta,\kappa)\leq \psi\bra{\frac{1}{2},\kappa}$ for all $\vartheta \in [0,1/2]$. The general case can be derived by substituting $T$ with
$\begin{pmatrix}
 T &0 \\
  0 & S \\
\end{pmatrix}$ and $X$ by $\begin{pmatrix}
 X &0 \\
  0 & 0 \\
\end{pmatrix}$.
\end{proof}
%=======================Remark2.4==================================
\begin{remark} By setting $\kappa$ to be equal to half (i.e., $\kappa=\frac{1}{2}$) in Theorem \ref{Rahma1}, we can deduce that we arrive at Theorem 2.3 as presented in \cite{KS}. Consequently, our findings represent an enhancement of the results established in that theorem.
\end{remark}
%=================================================================
As a consequence of Theorem \ref{Rahma1}, we have
%========================Corollary2.5=================================================================
\begin{corollary}\label{Rahma2} Let $T,S,X\in M_n(\c)$ with $T,S$ positive definite. Then for any unitarily invariant norm $\vertiii{\cdot}$ and a matrix monotone increasing function $\psi:(0,\infty)\to(0,\infty)$  with $\psi^*(x)=x(\psi(x))^{-1}$,
\begin{eqnarray*}
  && \frac{1}{2}\vertiii{T^{\frac{\mu}{2}}\bra{\psi(T^{\mu})X\psi^*(S^{\mu})+\psi^*(T^{\mu})X\psi(S^{\mu})}S^{\frac{\mu}{2}}}\\
 && \leq \vertiii{(1-\vartheta)T^{\kappa}XS^{1-\kappa}+\vartheta\bra{\frac{TX+XS}{2}}}.
\end{eqnarray*}
\end{corollary}
%==========================Corollary2.6=============================================
\begin{corollary}\label{cor2.5} Let $T,S,X\in M_n(\c)$ with $T,S$ positive definite. Then for any unitarily invariant norm $\vertiii{\cdot}$,
$ \frac{1}{4}\leq \mu\leq \frac{3}{4}$, $\kappa\in [0,1]$ and $\vartheta\in [1/2,\infty)$,
\begin{equation*}
  \frac{1}{2}\vertiii{T^{\mu}XS^{1-\mu}+T^{1-\mu}XS^{\mu}}\leq \vertiii{(1-\vartheta)T^{\kappa}XS^{1-\kappa}+\vartheta\bra{\frac{TX+XS}{2}}}.
\end{equation*}
\end{corollary}
\begin{proof} Letting $\psi(x)=\sqrt{x}$ in Corollary \ref{Rahma2}, we derived the result.
\end{proof}
%==============================================================
The following result is a consequence of Theorem \ref{Rahma1}.
%=======================corollary2.7================================================
\begin{corollary}\label{Boshra1} Let $T,S,X\in M_n(\c)$ with $T,S$ positive definite, $\eta=\min\{sp(T),sp(S)\}$,
$\mu\in [1/4,3/4]$ and $\kappa\in [0,1]$. Then for any unitarily invariant norm $\vertiii{\cdot}$
and a matrix monotone increasing function $\psi:(0,\infty)\to(0,\infty)$
\begin{eqnarray*}
 && \frac{\eta}{2f(\eta)}\vertiii{T^{\frac{\mu}{2}}\bra{\psi(T^{\mu})X+X\psi(S^{\mu})}S^{\frac{\mu}{2}}}\\
 &&\leq \vertiii{(1-\vartheta)T^{\kappa}XS^{1-\kappa}+\vartheta\bra{\frac{TX+XS}{2}}}
\end{eqnarray*}
holds for every $\vartheta\in [1/2,\infty)$.
\end{corollary}
%=======================================================================
Choosing $\psi(x)=\log(1+x)$ in Corollary \ref{Boshra1}, we have
%============================Corollary2.8==============================
\begin{corollary}
 Let $T,S,X\in M_n(\c)$ with $T,S$ positive definite, $\eta=\min\{sp(T),sp(S)\}$, $\mu\in [1/4,3/4]$ and $\kappa\in [0,1]$. Then for any unitarily invariant norm $\vertiii{\cdot}$
\begin{eqnarray*}
 && \frac{\eta}{2\log(1+\eta)}\vertiii{T^{\frac{\mu}{2}}\bra{\log(1+T^{\mu})X+X\log(1+S^{\mu})}S^{\frac{\mu}{2}}}\\
 &&\leq \vertiii{(1-\vartheta)T^{\kappa}XS^{1-\kappa}+\vartheta\bra{\frac{TX+XS}{2}}}
\end{eqnarray*}
holds for every $\vartheta\in [1/2,\infty)$.
\end{corollary}
%================================Theorem2.9======================================================
\begin{theorem}\label{Refund1} Consider $T,S,X\in M_n(\c)$ with $T$ and $S$ being positive definite, $\kappa\in [0,1]$, and $\vertiii{.}$ denoting any unitarily invariant norm. The function can be expressed as:
$$\phi(\vartheta,\kappa)=\vertiii{\bra{1-\frac{\vartheta}{2}}\bra{T^{\kappa}XS^{1-\kappa}+T^{1-\kappa}XS^{\kappa}}+\vartheta\bra{\frac{TX+XS}{2}}}$$
is increasing for $\frac{1}{2}\leq\vartheta<\infty$ and $\phi(\vartheta,\kappa)\leq \phi\bra{\frac{1}{2},\kappa}$
for all $\vartheta\in \sbra{0,\frac{1}{2}}$.
\end{theorem}
\begin{proof} Once again following the same lines of the proof of Theorem (\ref{Rahma1}), we
shall prove the result for $\vartheta>0$, $T=S$ and $T=diag(\eta_1,\cdots,\eta_n)$.
Suppose
$$\phi(\vartheta,\kappa)=\vertiii{\bra{1-\frac{\vartheta}{2}}\bra{T^{\kappa}XT^{1-\kappa}
+T^{1-\kappa}XT^{\kappa}}+\vartheta\bra{\frac{TX+XT}{2}}}=\frac{\vartheta}{2}Z(\vartheta,\kappa),$$
where $Z(\vartheta,\kappa)=\vertiii{W_1(\vartheta)\bra{T^{\kappa}XT^{1-\kappa}+T^{1-\kappa}XT^{\kappa}}+(TX+XT)}$ and $W_1(\vartheta)=\frac{2}{\vartheta}-1$.
\begin{eqnarray*}
  &&W_1(\vartheta)\bra{T^{\kappa}XT^{1-\kappa}+T^{1-\kappa}XT^{\kappa}}+(TX+XT)\\
&&=\sbra{\bra{W_1(\vartheta)\bra{\eta_i^{\mu}\eta_j^{1-\mu}+\eta_i^{1-\mu}\eta_j^{\mu}}+\eta_i+\eta_j}x_{ij}}_{i,j} \\
  &&=Y\circ \bra{W_1(\varrho)\bra{T^{\kappa}XT^{1-\kappa}+T^{1-\kappa}XT^{\kappa}}+TX+XS}.
\end{eqnarray*}
Now the matrix $Y$ can be written as
\begin{eqnarray*}
  &&\bra{\frac{W_1(\vartheta)\bra{\eta_i^{\mu}\eta_j^{1-\mu}+\eta_i^{1-\mu}\eta_j^{\mu}}+\eta_i+\eta_j}
{W_1(\varrho)\bra{\eta_i^{\mu}\eta_j^{1-\mu}+\eta_i^{1-\mu}\eta_j^{\mu}}+\eta_i+\eta_j}}_{i,j}
  =\bra{1+\frac{(W_1(\vartheta)-W_1(\varrho))\eta_i^{\kappa}\eta_j^{\kappa}}
{(W_1(\varrho)-1)\eta_i^{\kappa}\eta_j^{\kappa}+\eta_i^{1-\kappa}+\eta_j^{1-\kappa}}}_{i,j}\\
&&=(1)_{i,j}+\bra{\eta_i^{\kappa}\bra{\frac{W_1(\vartheta)-W_1(\varrho)}{(W_1(\varrho)-1)\eta_i^{\kappa}\eta_j^{\kappa}
+\eta_i^{1-\kappa}+\eta_j^{1-\kappa}}}\eta_j^{\mu}}_{i,j}
\end{eqnarray*}
Once again, considering Lemma (\ref{HHM2}), we observe that the latter matrix is positive semidefinite if and only if
 $W_1(\vartheta)>W_1(\varrho)$ and $2> W_1(\varrho)-1> -2$. Since $W(\vartheta)= W_1(\vartheta)-1= (2\vartheta^{-1}-2)$, it is a continuously decreasing function in the positive half-line and maps to the interval $\bras{-2,2}$ for $\vartheta$ in the range $[1/2,\infty)$. Therefore, as demonstrated in Theorem (\ref{Rahma1}), we can deduce that $W(\vartheta)>W(\varrho)$ and consequently, $W_1(\vartheta)>W_1(\varrho)$ for all $\vartheta\leq \varrho$.

Applying Theorem (\ref{HHM1}), we establish that $T(\vartheta,\kappa)\leq \frac{W_1(\vartheta)+1}{W_1(\varrho)+1}T(\varrho,\kappa)$. This verifies the result for the case when $T=S$ and $\vartheta\in [1/2,\infty)$.

For $\vartheta\in (0,1/2]$, we can observe that $3=W_1(1/2)\leq W_1(\vartheta)<\infty$, and according to Lemma (\ref{HHM2}), the matrix $Y$ with $\varrho=1/2$ is positive semidefinite. Similarly, the case $\vartheta=0$ can be established through the positive semidefiniteness of the matrix $\bra{\eta_i^{\kappa}\bra{\frac{W_1(\vartheta)-W_1(\varrho)}{(W_1(\varrho)-1)\eta_i^{\kappa}\eta_j^{\kappa}
+\eta_i^{1-\kappa}+\eta_j^{1-\kappa}}}\eta_j^{\mu}}_{i,j}$, which is confirmed by utilizing Lemma (\ref{HHM2}). This leads us to the desired result for this case, i.e., $\vartheta Z(\vartheta,\kappa)\leq \frac{1}{2}Z(\frac{1}{2},\kappa)$. In other words, $\phi(\vartheta,\kappa)\leq \phi(1/2,\kappa)$ holds for all $\vartheta \in[0,1/2]$.

The general case can be obtained by substituting $T$ with
$\begin{pmatrix}
 T &0 \\
  0 & S \\
\end{pmatrix}$ and $X$ by $\begin{pmatrix}
 X &0 \\
  0 & 0 \\
\end{pmatrix}$.
\end{proof}
%==========================================================================
The following outcome is an implication of Theorems \ref{Rahma1}, \ref{Refund1}, and Corollary \ref{cor2.5}, resulting in:
%==============================================================================
\begin{corollary} Let $T,S,X \in M_n(\c)$ with $T,S$ positive definite, $\kappa\in [0,1]$ and $\psi(\vartheta,\kappa)$ and
$\phi(\vartheta,\kappa)$ are same as taken in Theorem (\ref{Rahma1}) and (\ref{Refund1}) respectively. Then
\begin{equation}\label{Refund2}
  \psi(0,\kappa)\leq \frac{1}{2}\phi(0,\kappa)\leq \psi(\vartheta,\kappa)
\end{equation}
for $\vartheta\in \brass{1/2,\infty}$, or equivalently, for any unitarily invariant norm $\vertiii{\cdot}$
and $t\in \bras{-2,2}$,
\begin{equation*}
  \vertiii{T^{\kappa}XS^{1-\kappa}}\leq \frac{1}{2}\vertiii{T^{\kappa}XS^{1-\kappa}+T^{1-\kappa}XS^{\kappa}}\leq \frac{1}{t+2}\vertiii{TX+XS+tT^{\kappa}XS^{1-\kappa}}.
\end{equation*}
\end{corollary}
%=================================================================================
\begin{remark} (i) It's worth noting that the corollary mentioned earlier (\ref{Rahma2}) represents one of the potential enhancements to an inequality introduced by Kaur and Singh in their work (see \cite[Corollary 2.4]{KS}).\\
(ii) Take note that when we set $\kappa$ to the value of one-third (i.e., $\kappa=\frac{1}{3}$), it is evident that we arrive at the outcome outlined in Theorem 2.10 in \cite{KS}. This implies that our finding constitutes a broader and more generalized version of their result.
\end{remark}
%=======================================================================
\section{Sharpening of the Heinz inequalities and its reverses with the Kantorovich Constant}
%=============================================================================
%=====================================================================================
\indent In this section, we make a refinement of  Heinz inequality with the Kantorovich constant.
%=============================Lemma3.1=======================================
\begin{lemma}\label{med1}
  Let $\rho,\sigma>0$ and $0\leq \nu<\kappa\leq 1$. Then
\begin{equation}\label{B1}
  r(\sqrt{\rho\sharp_{\kappa}\sigma}-\sqrt{\sigma})^2+K(\sqrt{h},2)^{r'}\rho\sharp_{\nu}\sigma\leq
\nu \rho+(1-\nu)\sigma-\left(\frac{\nu}{\kappa}\right)(\rho\nabla_{\kappa}\sigma-\rho\sharp_{\kappa}\sigma),
\end{equation}
where $r=\min\{\frac{\nu}{\kappa},1-\frac{\nu}{\kappa}\}$, $h=\frac{\rho}{\sigma}$ and $r'=\min\{2r,1-2r\}$.
\end{lemma}
\begin{proof}
  An simple argument shows that
\begin{eqnarray}\label{B2}
  &&\nu \rho+(1-\nu)\sigma-\frac{\nu}{\kappa}\left(\rho\nabla_{\kappa}\sigma-\rho\sharp_{\kappa}\sigma\right)
=\nu \rho+(1-\nu)\sigma-\frac{\nu}{\kappa}\left(\kappa \rho+(1-\kappa)\sigma-\rho^{\kappa}\sigma^{1-\kappa}\right)\nonumber \\
 &&=\frac{\nu}{\kappa}\rho^{\kappa}\sigma^{1-\kappa}+\left(1-\frac{\nu}{\kappa}\right)\sigma=(\rho\sharp_{\kappa}\sigma)\nabla_{\frac{\nu}{\kappa}} \sigma.
\end{eqnarray}
By applying the inequality (\ref{A6}) for the relation (\ref{B2}), it follows that
$$r(\sqrt{\rho\sharp_{\kappa}\sigma}-\sqrt{\sigma})^2+K(\sqrt{h},2)^{r'}\rho\sharp_{\nu}\sigma\leq (\rho\sharp_{\kappa}\sigma)\nabla_{\frac{\nu}{\kappa}} \sigma. $$
Hence, the inequality (\ref{B1}) follows.
\end{proof}
%%%%%%%%%%%%%%%%=============================Lemma3.2================================
\begin{lemma}\label{med2}
   Let $\rho,\sigma>0$ and $0\leq \nu<\kappa\leq 1$. Then
\begin{equation}\label{B3}
   \nu \rho+(1-\nu)\sigma-\left(\frac{\nu}{\kappa}\right)(\rho\nabla_{\kappa}\sigma-\rho\sharp_{\kappa}\sigma)\leq K(\sqrt{h},2)^{-r'}\rho\sharp_{\nu}\sigma+R(\sqrt{\rho\sharp_{\kappa}\sigma}-\sqrt{\sigma})^2
\end{equation}
where $R=\max\{\frac{\nu}{\kappa},1-\frac{\nu}{\kappa}\}$, $h=\frac{\rho}{\sigma}$ and $r'=\min\{2r,1-2r\}$.
\end{lemma}
\begin{proof}
  By applying the inequality (\ref{A7}) for the relation (\ref{B2}), it follows that
$$\nu \rho+(1-\nu)\sigma-\left(\frac{\nu}{\kappa}\right)(\rho\nabla_{\kappa}\sigma-\rho\sharp_{\kappa}\sigma)\leq K(\sqrt{h},2)^{-r'}\rho\sharp_{\nu}\sigma+R(\sqrt{\rho\sharp_{\kappa}\sigma}-\sqrt{\sigma})^2.$$
So, we get the inequality (\ref{B3}).
\end{proof}
%======================================================================
\indent For two non-negative real numbers $\rho$ and $\sigma$, we define the Heinz mean in the parameter
$\mu$, $0\leq \mu\leq 1$,  as
\begin{equation}\label{Heinz1}
  H_{\mu}=\frac{\rho^{\mu}\sigma^{1-\mu}+\rho^{1-\mu}\sigma^{\mu}}{2}.
\end{equation}
Note that $H_{0}(\rho,\sigma)=H_{1}(\rho,\sigma)=\frac{\rho+\sigma}{2}$ and $H_{\frac{1}{2}}(\rho,\sigma)=\sqrt{\rho\sigma}$.  It is easy to
see that as a function of $\mu$, $H_{\mu}(\rho,\sigma)$ is convex, attains its minimum at $\mu=\frac{1}{2}$,
and attains its maximum at $\mu=0$ and $\mu=1$. Moreover, $H_{\mu}(\rho,\sigma)=H_{1-\mu}(\rho,\sigma)$ for
$0\leq \mu\leq 1$. Thus, the Heinz mean interpolates between the geometric mean and the arithmetic
mean:
\begin{equation}\label{Heinz2}
  \sqrt{\rho\sigma}\leq H_{\mu}(\rho,\sigma)\leq \frac{\rho+\sigma}{2}\,\,\, \mbox{for}\,\, 0\leq \mu\leq 1.
\end{equation}
%========================Theorem3.3=================================
\begin{theorem}\label{MUNA1} Let $\rho,\sigma>0$ and $0\leq \kappa<\kappa\leq 1$. Then
  \begin{equation}\label{Heinz-A1}
  \left.
    \begin{array}{ll}
     r\sbra{H_{\kappa}(\rho,\sigma)+H_0(\rho,\sigma)-H_{\frac{\kappa}{2}}(\rho,\sigma)}+K\sbra{\sqrt{h},2}^{r'}H_{\nu}(\rho,\sigma)\\
    \leq H_0(\rho,\sigma)-\bra{\frac{\nu}{\kappa}}\sbra{H_0(\rho,\sigma)-H_{\kappa}(\rho,\sigma)},
    \end{array}
  \right.
  \end{equation}
  where $r=\min\{\frac{\nu}{\kappa},1-\frac{\nu}{\kappa}\}$, $h=\frac{\rho}{\sigma}$ and $r'=\min\{2r,1-2r\}$.
\end{theorem}
\begin{proof} Interchanging $\rho$ with $\sigma$ and $\sigma$ with $\rho$ in inequality (\ref{B1}), we get
\begin{equation}\label{HF1}
  r(\sqrt{\sigma\sharp_{\kappa}\rho}-\sqrt{\rho})^2+K(\sqrt{h},2)^{r'}\sigma\sharp_{\nu}\rho\leq \nu \sigma+(1-\nu)\rho-\left(\frac{\nu}{\kappa}\right)(\sigma\nabla_{\kappa}\rho-\sigma\sharp_{\kappa}\rho).
\end{equation}
Adding (\ref{B1}) and (\ref{HF1}), we have
\begin{equation*}
  \left.
    \begin{array}{ll}
     r\sbra{(\sqrt{\rho\sharp_{\kappa}\sigma}-\sqrt{\sigma})^2+(\sqrt{\sigma\sharp_{\kappa}\rho}
-\sqrt{\rho})^2}+K\sbra{\sqrt{h},2}^{r'}\bra{2H_{\nu}(\rho,\sigma)}\\
\leq 2H_0(\rho,\sigma)-\bra{\frac{\nu}{\kappa}}\sbra{2H_0(\rho,\sigma)-2H_{\kappa}(\rho,\sigma)}
    \end{array}
  \right.
\end{equation*}
and so
\begin{equation*}
  \left.
    \begin{array}{ll}
     r\sbra{H_{\kappa}(\rho,\sigma)+H_0(\rho,\sigma)-H_{\frac{\kappa}{2}}(\rho,\sigma)}+K\sbra{\sqrt{h},2}^{r'}H_{\nu}(\rho,\sigma)\\
\leq H_0(\rho,\sigma)-\bra{\frac{\nu}{\kappa}}\sbra{H_0(\rho,\sigma)-H_{\kappa}(\rho,\sigma)}.
    \end{array}
  \right.
\end{equation*}
\end{proof}
In similar of proof of Theorem \ref{MUNA1}, we can prove the following result.
%===============================Theorem3.4===================================
\begin{theorem}
  Let $\rho,\sigma>0$ and $0\leq \kappa<\kappa\leq 1$. Then
  \begin{equation}\label{Heinz-A2}
  \left.
    \begin{array}{ll}
    H_0(\rho,\sigma)-\bra{\frac{\nu}{\kappa}}\sbra{H_0(\rho,\sigma)-H_{\kappa}(\rho,\sigma)}\leq K\sbra{\sqrt{h},2}^{-r'}H_{\nu}(\rho,\sigma)\\
+R\sbra{H_{\kappa}(\rho,\sigma)+H_0(\rho,\sigma)-H_{\frac{\kappa}{2}}(\rho,\sigma)},
    \end{array}
  \right.
  \end{equation}
  where $R=\min\{\frac{\nu}{\kappa},1-\frac{\nu}{\kappa}\}$, $h=\frac{\rho}{\sigma}$ and $r'=\min\{2r,1-2r\}$.
\end{theorem}
%===================section4=============================================
\section{New operator versions of Heinz-type inequalities}
%====================================================================================
Let $\h$ represent a complex Hilbert space, and $\bh$ denote the $C^*$-algebra comprising all bounded linear operators on $\h$. An operator $T\in\bh$ is considered positive if $\seq{Tx,x}\geq 0$ holds true for every $x\in\h$. We express this as $T\geq 0$. \
Now, let $T$ and $S$ be two positive operators in $\bh$, and $\kappa$ take on values in the interval $[0,1]$. The $\kappa$-weighted arithmetic mean of $T$ and $S$, denoted as $T\nabla_{\kappa} S$, is defined as:
\begin{equation*}
T\nabla_{\kappa} S=(1-\kappa)T+\kappa S.
\end{equation*}
When $T$ is invertible, the $\kappa$-geometric mean of $T$ and $S$, represented as $T\sharp_{\kappa}S$, is defined as:
\begin{equation*}
T\sharp_{\kappa}S=T^{\frac{1}{2}}\left(T^{-\frac{1}{2}}ST^{-\frac{1}{2}}\right)^{\kappa}T^{\frac{1}{2}}.
\end{equation*}
In the case where $\kappa=\frac{1}{2}$, we can simplify the notation to $T\nabla S$ and $T\sharp S$ to refer to the $\kappa$-weighted arithmetic mean and the $\kappa$-geometric mean, respectively. It is well-known that for positive invertible operators $T$ and $S$, the following inequality holds:
\begin{equation*}
T\sharp_{\kappa}S\leq T \nabla_{\kappa}S, \quad \kappa\in [0,1].
\end{equation*}
Additionally, we define the operator version of the Heinz mean as $H_{\kappa}(T,S)$:
\begin{equation*}
H_{\kappa}(T,S)=\frac{T\sharp_{\kappa}S+T\sharp_{1-\kappa}S}{2}
\end{equation*}
for the case where $T$ and $S$ are positive invertible operators and $\kappa\in[0,1]$.\\
%==============================================================================
\indent In this section, we will present improved variants of Heinz-type operator inequalities and their converses, exploiting the monotonicity of operator functions as the foundational concept for the ensuing discussion.
%================================lemma4.1==========================================
\begin{lemma}\label{lemma4.1}\cite{PFMHS} Suppose $T\in\bh$ is self-adjoint. If $f$ and $g$ are continuous functions such that $f (t)\geq g(t)$ for $t\in sp(T)$ (where $sp(T)$ represents the spectrum of the operator $T$), then it follows that $f (T)\geq g(T)$.
\end{lemma}
%=======================================================================================
Next we present our main results on the basis of inequality (\ref{Heinz-A1}). By Lemma \ref{lemma4.1},
we have the following.
%=========================Theorem4.2=====================================
\begin{theorem}\label{thm4.2} Let $T,S\in\bh$ be positive invertible operators, $I$ is the identity operator and
$0\leq \kappa<\kappa\leq 1$. If all positive numbers $m, m'$ and $M, M'$ satisfy either of the conditions
$0<mI\leq T\leq m'I<M'I\leq S\leq MI$  or  $0< mI \leq  S \leq m' I\leq T\leq MI$, then:
\begin{equation}\label{Heinz-O1}
  \left.
    \begin{array}{ll}
     r\sbra{H_{\kappa}(T,S)+H_0(T,S)-H_{\frac{\kappa}{2}}(T,S)}+K\sbra{\sqrt{h},2}^{r'}H_{\nu}(T,S)\\
    \leq H_0(T,S)-\bra{\frac{\nu}{\kappa}}\sbra{H_0(T,S)-H_{\kappa}(T,S)},
    \end{array}
  \right.
  \end{equation}
  where $r=\min\{\frac{\nu}{\kappa},1-\frac{\nu}{\kappa}\}$, $h=\frac{M}{m}$ and $r'=\min\{2r,1-2r\}$.
\end{theorem}
\begin{proof} Assuming that $0 \leq \nu < \kappa \leq 1$, according to inequality (\ref{Heinz-A1}), for any positive value of $x$, we can conclude:
\begin{equation*}
  \left.
    \begin{array}{ll}
     r\sbra{H_{\kappa}(1,x)+H_0(1,x)-H_{\frac{\kappa}{2}}(1,x)}+K\sbra{\sqrt{h},2}^{r'}H_{\nu}(1,x)\\
    \leq H_0(1,x)-\bra{\frac{\nu}{\kappa}}\sbra{H_0(1,x)-H_{\kappa}(1,x)},
    \end{array}
  \right.
  \end{equation*}
Regarding the operator $X = T^{-1/2}ST^{-1/2}$, within the framework of the first condition, we establish the following range: $I \leq hI = \frac{M}{m}I \leq X \leq h' I = \frac{M'}{m'}I$. Consequently, we infer that $\sigma(X) \subseteq [h,h'] \subseteq (1,\infty)$. Applying Lemma \ref{lemma4.1}, we obtain:
\begin{eqnarray*}
 &&r\sbra{H_{\kappa}(I,X)+H_0(I,X)-H_{\frac{\kappa}{2}}(I,X)}+\min_{h\leq x\leq h'} K\sbra{\sqrt{x},2}^{r'}H_{\nu}(I,X) \\
  &&\leq H_0(I,X)-\bra{\frac{\nu}{\kappa}}\sbra{H_0(I,X)-H_{\kappa}(I,X)},
\end{eqnarray*}
As the Kantorovich constant $K(t,2) = \frac{(1+t)^2}{4t}$ exhibits monotonicity within the interval $(0,\infty)$, it follows that:
\begin{eqnarray}\label{O2}
 &&r\sbra{H_{\kappa}(I,T^{-1/2}ST^{-1/2})+H_0(I,T^{-1/2}ST^{-1/2})-H_{\frac{\kappa}{2}}(I,T^{-1/2}ST^{-1/2})}\nonumber \\
  && +\min_{h\leq x\leq h'} K\sbra{\sqrt{x},2}^{r'}H_{\nu}(I,T^{-1/2}ST^{-1/2})\leq H_0(I,T^{-1/2}ST^{-1/2})\nonumber\\
&&-\bra{\frac{\nu}{\kappa}}\sbra{H_0(I,T^{-1/2}ST^{-1/2})-H_{\kappa}(I,T^{-1/2}ST^{-1/2})},
\end{eqnarray}
Likewise, within the context of the second condition, we observe that $I \leq \frac{1}{h} I = \frac{m}{M} h \leq X \leq \frac{1}{h'} I = \frac{m'}{M'} I$. Utilizing Lemma \ref{lemma4.1}, we obtain the following:
\begin{equation*}
  \left.
    \begin{array}{ll}
     r\sbra{H_{\kappa}(I,X)+H_0(I,X)-H_{\frac{\kappa}{2}}(I,X)}+\min_{\frac{1}{h'}\leq x\leq \frac{1}{h}} K\sbra{\sqrt{x},2}^{r'}H_{\nu}(I,X)\\
    \leq H_0(I,X)-\bra{\frac{\nu}{\kappa}}\sbra{H_0(I,X)-H_{\kappa}(I,X)},
    \end{array}
  \right.
  \end{equation*}
Since the Kantorovich constant $K(t,2)=\frac{(1+t)^2}{4t}$ is an increasing function on $(0,\infty)$, then
\begin{eqnarray*}\label{O3}
  && r\sbra{H_{\kappa}(I,T^{-1/2}ST^{-1/2})+H_0(I,T^{-1/2}ST^{-1/2})-H_{\frac{\kappa}{2}}(I,T^{-1/2}ST^{-1/2})}\nonumber \\
  && +\min_{\frac{1}{h'}\leq x\leq \frac{1}{h}} K\sbra{\sqrt{x},2}^{r'}H_{\nu}(I,T^{-1/2}ST^{-1/2})\leq H_0(I,T^{-1/2}ST^{-1/2})\nonumber\\
&&-\bra{\frac{\nu}{\kappa}}\sbra{H_0(I,T^{-1/2}ST^{-1/2})-H_{\kappa}(I,T^{-1/2}ST^{-1/2})},
\end{eqnarray*}
By multiplying both inequalities (\ref{O2}) and (\ref{O3}) on both the left-hand and right-hand sides by the operator $T^{1/2}$, we can infer the desired inequality (\ref{Heinz-O1}).
\end{proof}
%=============================================================================
\begin{theorem}\label{thm4.3} Consider positive invertible operators $T$ and $S$ in a Hilbert space $\mathcal{H}$, where $I$ represents the identity operator. Additionally, let $\kappa$ be a non-negative number such that $0 \leq \kappa < \kappa \leq 1$. Assuming that there exist positive real numbers $m, m', M, M'$ that satisfy either of the following conditions:
\begin{enumerate}
  \item [(a)] $0 < mI \leq T \leq m'I < M'I \leq S \leq MI$
  \item [(b)] $0 < mI \leq S \leq m'I \leq T \leq MI$
\end{enumerate}
Then, the following conclusions hold:
\begin{equation}\label{Heinz-A2}
  \left.
    \begin{array}{ll}
    H_0(T,S)-\bra{\frac{\nu}{\kappa}}\sbra{H_0(T,S)-H_{\kappa}(T,S)}\leq K\sbra{\sqrt{h},2}^{-r'}H_{\nu}(T,S)\\
+R\sbra{H_{\kappa}(T,S)+H_0(T,S)-H_{\frac{\kappa}{2}}(T,S)},
    \end{array}
  \right.
  \end{equation}
  where $R=\min\{\frac{\nu}{\kappa},1-\frac{\nu}{\kappa}\}$, $h=\frac{M}{m}$ and $r'=\min\{2r,1-2r\}$.
\end{theorem}
\begin{proof}The proof process is similar to that of Theorem \ref{thm4.2}, and thus, we will not provide it here.
\end{proof}
%===============================================================================
\begin{remark} The nature of the Kantorovich constant's characteristics makes it clear that the inequalities outlined in Theorems \ref{thm4.2} and \ref{thm4.3} signify improved results compared to those detailed in \cite{Rash1}, \cite{RashNa}, \cite{SS}, \cite{YL}, and \cite{ZW}.
\end{remark}
%==============================================================
\section{Utilizations of the improved Young-type inequalities for traces, determinants, and norms of positive definite matrices}
%========================================================================
In this section, we introduce a collection of improved Young-type inequalities designed specifically for traces, determinants, and norms of positive semi-definite matrices.\\
%=========================================================
\indent A matrix version proved in \cite{Ando} says that if $T, S \in M_n(\c)$ are positive semi-definite, then
\begin{equation}\label{Ando1}
  s_j(TS)\leq s_j\bra{\frac{1}{p}T^p+\frac{1}{q}S^q}
\end{equation}
for $j=1,\cdots,n$
%=========================lemma5.1=======================================================================
\begin{lemma}
  Let $\rho,\sigma>0$, $p,q>1$ such that $\frac{1}{p}+\frac{1}{q}=1$. Then for $m\in\N$, we have
  \begin{equation}\label{Man1}
    \bra{\rho^{\frac{1}{p}}\sigma^{\frac{1}{q}}}^{m}+r_0^m\bra{\rho^{\frac{m}{2}}-\sigma^{\frac{m}{2}}}^2\leq \bra{\frac{\rho^r}{p}+\frac{\sigma^r}{q}}^{\frac{m}{r}},\quad r\geq1
  \end{equation}
  where $r_0=\min\{\frac{1}{p},\frac{1}{q}\}$.
\end{lemma}
%===========================lemma5.2=============================
\begin{lemma}\label{Benz1}
  Let $T_i\in M_n(\c)$ $(i=1,\cdots,n)$,. Then
  $$\sum_{j=1}^{n}s_j(T_1\cdots T_n)\leq \sum_{j=1}^{n}s_j(T_1)\cdots s_j(T_k).$$
\end{lemma}
%==========================THEOREM5.3===============================================
\begin{theorem}\label{TRACEONE}
  Let $T,S\in\bh$ be positive definite, $p,q>1$ such that $\frac{1}{p}+\frac{1}{q}=1$ and $m\in\N$.
  Then
  \begin{equation}\label{Man2}
    \bra{\frac{tr(T^r)}{p}+\frac{tr(S^r)}{q}}^{\frac{m}{r}}\geq
    \bra{tr\abs{T^{\frac{1}{p}}S^{\frac{1}{q}}}}^m+r_0^m\bra{\bra{tr(T)}^{\frac{m}{2}}-\bra{tr(S)}^{\frac{m}{2}}}^{2},
  \end{equation}
  where $r_0=\min\{\frac{1}{p},\frac{1}{q}\}$.
\end{theorem}
\begin{proof}
  By inequality (\ref{Man2}), we have
  \begin{eqnarray*}
    s_j^{\frac{m}{r}}\bra{\frac{T^r}{p}+\frac{S^r}{q}} &=& \bra{\frac{s_j^{\frac{m}{r}}(T^r)}{p}+\frac{s_j^{\frac{m}{r}}(S^r)}{q}} \\
     &\geq& s_j^{m}\bra{T^{\frac{1}{p}}} s_j^m\bra{S^{\frac{1}{q}}}+r_0^m\bra{s_j^{\frac{m}{2}}(T)-s_j^{\frac{m}{2}}(S)}^2\\
     &=& s_j^{m}\bra{T^{\frac{1}{p}}} s_j^m\bra{S^{\frac{1}{q}}}+r_0^m\bra{s_j^m(T)+s_j^m(S)-2s_j^{\frac{m}{2}}(T)s_j^{\frac{m}{2}}(S)}
  \end{eqnarray*}
  for $j =1,\cdots,n$. Thus, by Lemma \ref{Benz1} and the  Cauchy-Schwarz inequality, we have
  \begin{eqnarray*}
  tr^{\frac{m}{r}}\bra{\frac{T^r}{p}+\frac{S^r}{q}}  &=&\sum_{j=1}^{n}s_j^{\frac{m}{r}}\bra{\frac{T^r}{p}+\frac{S^r}{q}}\\
    &\geq&\sum_{j=1}^{n}s_j^{m}\bra{T^{\frac{1}{p}}}s_j^m\bra{S^{\frac{1}{q}}} \\
    &+&r_0^m\bra{\sum_{j=1}^{n}s_j^m(T)+\sum_{j=1}^{n}s_j^m(S)-2\sum_{j=1}^{n}s_j^{\frac{m}{2}}(T)s_j^{\frac{m}{2}}(S)} \\
    \end{eqnarray*}
    Hence
    \begin{eqnarray*}
    tr^{\frac{m}{r}}\bra{\frac{T^r}{p}+\frac{S^r}{q}}  &\geq&\sum_{j=1}^{n}s_j^{m}\bra{T^{\frac{1}{p}}S^{\frac{1}{q}}}\\
     &+&r_0^m\bra{\sum_{j=1}^{n}s_j^m(T)+\sum_{j=1}^{n}s_j^m(S)-2\sum_{j=1}^{n}s_j^{\frac{m}{2}}(T)s_j^{\frac{m}{2}}(S)} \\
     &\geq& \bra{tr\abs{\bra{T^{\frac{1}{p}}S^{\frac{1}{q}}}}}^{m}+r_0^m\left[(tr(T))^m+(tr(S))^m\right.\\
     &-&\left.2\bra{\sum_{j=1}^{n}s_j(T)}^{\frac{m}{2}}
     \bra{\sum_{j=1}^{n}s_j(S)}^{\frac{m}{2}}\right]\\
     &=&\bra{tr\abs{\bra{T^{\frac{1}{p}}S^{\frac{1}{q}}}}}^{m}+r_0^m\bra{\bra{tr(T)}^{\frac{m}{2}}-\bra{tr(S)}^{\frac{m}{2}}}^2
  \end{eqnarray*}
\end{proof}
%==========================remark5.4================================================
\begin{remark} Ando's singular value inequality (\ref{Ando1}) entails the norm inequality
\begin{equation}\label{Re1}
  \vertiii{T^{\kappa}S^{1-\kappa}}\leq \vertiii{\kappa T+(1-\kappa)S}.
\end{equation}
  So, our Theorem \ref{TRACEONE} improves this inequality for the trace norm:
  \begin{equation}\label{RASHID-1}
    \norm{T^{\frac{1}{p}}S^{\frac{1}{q}}}_1^m+r_0^m\bra{\norm{T}_1^{\frac{m}{2}}-\norm{S}_1^{\frac{m}{2}}}^2
  \leq \norm{\frac{1}{p}T^{r}+\frac{1}{q}S^{r}}_1^{\frac{m}{r}}
  \end{equation}
\end{remark}
%==============================THEOREM5.5=============================================
\begin{theorem}
  Let $T,S\in\bh$ be positive definite, $p,q>1$ such that $\frac{1}{p}+\frac{1}{q}=1$ and $m\in\N$.
  Then for all $r\geq 1$
  \begin{equation}\label{RashidQ1}
    \det\bra{\frac{T^r}{p}+\frac{S^r}{q}}^{\frac{m}{r}}\geq \det\bra{T^{\frac{1}{p}}S^{\frac{1}{q}}}^{m}
    +r_0^{mn}\det\bra{T^m+S^m-2S^{\frac{m}{2}}\bra{S^{-\frac{1}{2}}TS^{-\frac{1}{2}}}^{m}S^{\frac{m}{2}}}^2,
  \end{equation}
  where $r_0=\min\{\frac{1}{p},\frac{1}{q}\}$.
\end{theorem}
\begin{proof}
  By inequality (\ref{Man1}), we have
  $$s_j^{\frac{m}{r}}\bra{\frac{1}{p}\bra{S^{-\frac{r}{2}}T^rS^{-\frac{r}{2}}}+\frac{1}{q}I}
  \geq s_j^{\frac{m}{p}}\bra{S^{-\frac{1}{2}}TS^{-\frac{1}{2}}}+r_0^{m}\bra{s_j^{\frac{m}{2}}\bra{S^{-\frac{1}{2}}TS^{-\frac{1}{2}}}-1}^{2}$$
  for all $j=1,\cdots,n$.
  \begin{eqnarray*}
    \det\bra{\frac{1}{p}S^{-\frac{r}{2}}TS^{-\frac{r}{2}}+\frac{1}{q}}^{\frac{m}{r}} &=& \prod_{j=1}^{n}
    \bra{\frac{1}{p}s_j^{\frac{m}{r}}\bra{S^{-\frac{r}{2}}T^rS^{-\frac{r}{2}}+\frac{1}{q}}} \\
    &\geq& \prod_{j=1}^{n}\sbra{s_{j}^{\frac{m}{p}}\bra{S^{-\frac{1}{2}}TS^{-\frac{1}{2}}}
    +r_0^{m}\bra{s_j^{\frac{m}{2}}\bra{S^{-\frac{1}{2}}TS^{-\frac{1}{2}}}-1}^{2}}\\
    &\geq& \prod_{j=1}^{n}\sbra{s_{j}^{\frac{1}{p}}\bra{S^{-\frac{1}{2}}TS^{-\frac{1}{2}}}^{m}}\\
    &+&r_0^{mn}\prod_{j=1}^{n}\sbra{s_j^{\frac{m}{2}}\bra{S^{-\frac{1}{2}}TS^{-\frac{1}{2}}}-1}^{2}\\
   &=&\det\bra{S^{-\frac{1}{2}}TS^{-\frac{1}{2}}}^{\frac{m}{p}}+r_0^{mn}\sbra{\bra{S^{-\frac{1}{2}}TS^{-\frac{1}{2}}}^{\frac{m}{2}}-I}^2.
  \end{eqnarray*}
  Consequently
  $$ \det\bra{\frac{T^r}{p}+\frac{S^r}{q}}^{\frac{m}{r}}\geq \det\bra{T^{\frac{1}{p}}S^{\frac{1}{q}}}^{m}
    +r_0^{mn}\det\bra{T^m+S^m-2S^{\frac{m}{2}}\bra{S^{-\frac{1}{2}}TS^{-\frac{1}{2}}}^{m}S^{\frac{m}{2}}}^2.$$
\end{proof}
%======================THEOREM5.6=======================================
\begin{theorem}\label{GHADEER11}
  Let $T,S,X\in M_n(\c)$ such that $T$ and $S$ are positive semi-definite
  and $p,q>1$ with $\frac{1}{p}+\frac{1}{q}=1$ and $m\in\N$. Then for all $r\geq 1$, we have
\begin{equation}\label{un-Norm}
    \vertiii{T^{\frac{1}{p}}XS^{\frac{1}{q}}}^m+r_0^m\bra{\vertiii{TX}^{\frac{m}{2}}-\vertiii{SX}^{\frac{m}{2}}}^2
    \leq \bra{\frac{1}{p}\vertiii{TX}^{r}+\frac{1}{q}\vertiii{XB}^{r}}^{\frac{m}{r}},
  \end{equation}
  where $r_0=\min\{\frac{1}{p},\frac{1}{q}\}$.
\end{theorem}
%=============================================================
To prove Theorem \ref{GHADEER11}, we need the following lemma which is known as the
Heinz-Kato type for unitarily invariant norm.
%====================================================Lemma5.7==============
\begin{lemma}[\cite{kit0}]\label{A00}
  Let $T,S\in M_n(\c)$ be positive definite matrices and $0\leq \vartheta\leq 1$. Then we have
  \begin{equation}\label{A11}
    \vertiii{T^{\vartheta}XS^{1-\vartheta}}\leq \vertiii{TX}^{\vartheta}\vertiii{XB}^{1-\vartheta}.
  \end{equation}
  In particular
  \begin{equation}\label{A22}
    tr\abs{T^{\vartheta}XS^{1-\vartheta}}\leq \bra{tr(T)}^{\vartheta}\bra{tr(S)}^{1-\vartheta}.
  \end{equation}
\end{lemma}
%======================proof Theorem5.6===========================
\begin{proof}[Proof of Theorem \ref{GHADEER11}]
We have
\begin{eqnarray*}
   \vertiii{T^{\frac{1}{p}}XS^{\frac{1}{q}}}^m&+&r_0^m\bra{\vertiii{TX}^{\frac{m}{2}}-\vertiii{SX}^{\frac{m}{2}}}^2\\
   &\leq&\sbra{\vertiii{TX}^{\frac{1}{p}}\vertiii{XB}^{\frac{1}{q}}}^m+r_0^m\bra{\vertiii{TX}^{\frac{m}{2}}-\vertiii{SX}^{\frac{m}{2}}}^2 \\
   && (\mbox{by Lemma \ref{A00}}) \\
   &\leq& \bra{\frac{1}{p}\vertiii{TX}^{r}+\frac{1}{q}\vertiii{XB}^{r}}^{\frac{m}{r}}\,\,\,(\mbox{by inequality \ref{Man1}}).
\end{eqnarray*}
\end{proof}
%================================LEmma5.8==================================================
\begin{lemma}[\cite{furu}]
  Let $\omega_1,\cdots,\omega_n$ be non-negative real numbers and $\vartheta_1,\cdots,\vartheta_n$ be positive real numbers
  with $\sum_{i=1}^{n}\vartheta_i=1$. Then we have
  \begin{equation}\label{Furu}
  \displaystyle  \prod_{k=1}^{n}\omega_{k}^{\vartheta_k}+r\bra{\sum_{k=1}^{n}\omega_k-n\sqrt[n]{\prod_{k=1}^{n}\omega_k}}\leq \sum_{i=1}^{n}\vartheta_k\omega_k,
  \end{equation}
  where $\displaystyle r=\min\set{\vartheta_k:k=1,\cdots,n}$.
\end{lemma}
%=======================THEOREM5.9=====================================================================
\begin{theorem}\label{Rashid-2}
  Let $T_i\in M_n(\c)$ $(i=1,\cdots,n)$ be positive semi-definite. If $0\leq \vartheta_i\leq 1$ $(i=1,\cdots,n)$ with $\sum_{i=1}^{n}\vartheta_i=1$,
   then
   \begin{equation}\label{Trace1}
  \sum_{k=1}^{n}tr(\vartheta_kT_k)\geq tr\abs{\prod_{k=1}^{n}T_k^{\vartheta_k}}+r\bra{\sum_{k=1}^{n}tr(T_k)-n\sqrt[n]{\prod_{k=1}^{n}tr(T_k)}},
   \end{equation}
    where $\displaystyle r=\min\set{\vartheta_k:k=1,\cdots,n}$.
\end{theorem}
\begin{proof}
  By inequality (\ref{Furu}), we have
  $$\displaystyle \sum_{k=1}^{n}\vartheta_ks_j(T_k)\geq   \prod_{k=1}^{n}s_j(T_k)^{\vartheta_k}+r\bra{\sum_{k=1}^{n}s_j(T_k)-n\sqrt[n]{\prod_{k=1}^{n}s_j(T_k)}}$$
  for $j =1,\cdots,n$. Thus, by Lemma \ref{Benz1} and the generalized Cauchy-Schwarz inequality, we have
  \begin{eqnarray*}
    tr\bra{\sum_{k=1}^{n}\vartheta_kT_k} &=&\sum_{k=1}^{n}\vartheta_ktr(T_k)= \sum_{k=1}^{n}\vartheta_k\sum_{j=1}^{n}s_j(T_k)=\sum_{j=1}^{n}\sum_{k=1}^{n}\vartheta_ks_j(T_k)\\
     &\geq& \sum_{j=1}^{n}s_j(T_1^{\vartheta_1})\cdots s_j(T_k^{\vartheta_n})\\
     &+& r\bra{\sum_{j=1}^{n}\sum_{k=1}^{n}s_j(T_k)-n\sum_{j=1}^{n}\sqrt[n]{\prod_{k=1}^{n}s_j(T_n)}} \\
      &\geq& \sum_{j=1}^{n}s_j(T_1^{\vartheta_1}\cdots T_n^{\vartheta_n})\\
     &+& r\bra{\sum_{j=1}^{n}\sum_{k=1}^{n}s_j(T_k)-n\sqrt[n]{\prod_{k=1}^{n}\sum_{j=1}^{n}s_j(T_k)}} \\
     &\geq& tr\abs{T_1^{\vartheta_1}\cdots T_n^{\vartheta_n}}+r\bra{\sum_{k=1}^{n}tr(T_k)-n\sqrt[n]{\prod_{k=1}^{n}tr(T_k)}}
  \end{eqnarray*}
  where $\displaystyle r=\min\set{\vartheta_k:k=1,\cdots,n}$.
\end{proof}
Our Theorem \ref{Rashid-2} entils the following trace norm
\begin{equation}\label{Tracethree}
  \norm{\sum_{k=1}^{n}\vartheta_kT_k}_{1}\geq \norm{\prod_{k=1}^{n}T_k^{\vartheta_k}}_1+r\bra{\norm{\sum_{k=1}^{n}T_k}_1
  -n\sqrt[n]{\prod_{k=1}^{n}\norm{T_k}_1}},
\end{equation}
where $\displaystyle r=\min\set{\vartheta_k:k=1,\cdots,n}$.
%=================Theorem5.10===================================================
\begin{theorem}\label{DETERMINE}
  Let $T_i\in M_n(\c)$ $(i=1,\cdots,n)$ be positive definite. If $0\leq \vartheta_i\leq 1$ $(i=1,\cdots,n)$ with $\sum_{i=1}^{n}\vartheta_i=1$,
   then
   \begin{equation}\label{DETER1}
  \det\bra{\sum_{k=1}^{n}\vartheta_kT_k}\geq \prod_{k=1}^{n}\det\bra{T_k^{\vartheta_k}}+r\bra{\det\bra{\sum_{k=1}^{n}T_k}-n\sqrt[n]{\prod_{k=1}^{n}\det(T_k)}},
   \end{equation}
    where $\displaystyle r=\min\set{\vartheta_k:k=1,\cdots,n}$.
\end{theorem}
To prove Theorem \ref{DETERMINE}, we need the following lemma.
%==================Lemma5.11================================
\begin{lemma}[\cite{Horn}]\label{nth-Deter}
  Let $T,S\in M_n(\c)$ be positive definite. Then we have
  \begin{equation}\label{sum-det}
    \det(T+S)^{\frac{1}{n}}\geq \det(T)^{\frac{1}{n}}+\det(S)^{\frac{1}{n}}.
  \end{equation}
\end{lemma}
%===================Proof of THEOREM5.10=============================
\begin{proof}[Proof of Theorem \ref{DETERMINE}] We have
\begin{eqnarray*}
   \det\bra{\sum_{k=1}^{n}\vartheta_kT_k} &=& \sbra{\det\bra{\sum_{k=1}^{n}\vartheta_kT_k}^{\frac{1}{n}}}^{n} \\
   &\geq&\sbra{\sum_{k=1}^{n}\det\bra{\vartheta_kT_k}^{\frac{1}{n}}}^{n}\quad(\mbox{by Lemma \ref{nth-Deter}} ) \\
   &\geq&\sbra{\sum_{k=1}^{n}\vartheta_k\det\bra{T_k}^{\frac{1}{n}}}^{n}  \\
   &\geq & \sbra{\prod_{k=1}^{n}\bra{\bra{T_k}^{\frac{1}{n}}}^{\vartheta_k}}^{n}+
   r^n\bra{\sum_{k=1}^{n}{\det\bra{T_k}^{\frac{1}{n}}}-n\sqrt[n]{\det{\prod_{k=1}^{n}T_k}}} \\
   &=& \prod_{k=1}^{n}\det\bra{T_k^{\vartheta_k}}+
   r^n\bra{\sum_{k=1}^{n}{\det\bra{T_k}^{\frac{1}{n}}}-n\sqrt[n]{\det{\prod_{k=1}^{n}T_k}}}
\end{eqnarray*}
\end{proof}
%================================================================================
%======================Lemma5.12========================================================
\begin{lemma}[\cite{PD}]
  Let $\gamma_1,\gamma_2,\cdots,\gamma_n$ be a set of non-negative real numbers constrained by $\displaystyle\sum_{k=1}^{j}\gamma_k=\Gamma_j$.
  If $\omega_1,\omega_2,\cdots,\omega_n$ are positive real numbers, then
  \begin{equation}\label{Constrained}
    \frac{1}{\Gamma_n}\sum_{k=1}^{n}\gamma_k\omega_k+\sqrt{1+\bra{\frac{1}{\Gamma_n}\sum_{k=1}^{n}\gamma_k\omega_k}^2}
    \geq \sbra{\prod_{k=1}^{n}\bra{\omega_k+\sqrt{1+\omega_k^2}}^{\gamma_k}}^{\frac{1}{\Gamma_n}}
  \end{equation}
  holds.
\end{lemma}
%==================THEOREM5.13====================================
\begin{theorem}
  Let $T_1,\cdots,T_k\in M_n(\c)$ be positive define and let $\gamma_1,\gamma_2,\cdots,\gamma_n$ be a set of non-negative real numbers
  such that $\displaystyle\sum_{k=1}^{n}\gamma_k=\Gamma_n$. Then
  \begin{equation}\label{Const1}
    \frac{1}{\Gamma_n}\sum_{k=1}^{n}tr(T_k)+\sqrt{1+\bra{\frac{1}{\Gamma_n}\sum_{k=1}^{n}tr(T_k)}^2}
    \geq \prod_{k=1}^{n}\sbra{tr(T_k)+\sqrt{1+tr(T_k^2)}}^{ \frac{1}{\Gamma_n}}.
  \end{equation}
\end{theorem}
\begin{proof}
  By inequality (\ref{Constrained}), we have
  \begin{equation}\label{const2}
    \frac{1}{\Gamma_n}\sum_{k=1}^{n}\gamma_ks_j(T_k)+\sqrt{1+\bra{\frac{1}{\Gamma_n}\sum_{k=1}^{n}\gamma_ks_j(T_k)}^2}
    \geq \sbra{\prod_{k=1}^{n}\bra{s_j(T_k)+\sqrt{1+s_j^2(T_k)}}^{\gamma_k}}^{\frac{1}{\Gamma_n}}
  \end{equation}
  for all $j=1,\cdots,n$. Hence we have
  \begin{eqnarray*}
   \frac{1}{\Gamma_n}\sum_{k=1}^{n}\gamma_k\sum_{j=1}^{n} s_j(T_k)&+&\sqrt{1+\bra{\frac{1}{\Gamma_n}\sum_{k=1}^{n}\gamma_k\sum_{j=1}^{n}  s_j(T_k)}^2} \\
     &\geq& \sbra{\prod_{k=1}^{n}\bra{\sum_{j=1}^{n}s_j(T_k)+\sqrt{1+\sum_{j=1}^{n}s_j^2(T_k)}}^{\gamma_k}}^{\frac{1}{\Gamma_n}}\\
  \end{eqnarray*}
  Consequently,
  \begin{eqnarray*}
    tr\bra{\frac{1}{\Gamma_n}\sum_{k=1}^{n}\gamma_kT_k}&+&\sqrt{1+\bra{\frac{1}{\Gamma_n}\sum_{k=1}^{n}\gamma_k tr(T_k)}^2}\\
    &=&\frac{1}{\Gamma_n}\sum_{k=1}^{n}\gamma_ktr(T_k)+\sqrt{1+\bra{\frac{1}{\Gamma_n}\sum_{k=1}^{n}\gamma_k tr(T_k)}^2}\\
    &\geq& \sbra{\prod_{k=1}^{n}\bra{tr(T_k)+\sqrt{1+\sum_{j=1}^{n}tr(T_k^2)}}^{\gamma_k}}^{\frac{1}{\Gamma_n}}
  \end{eqnarray*}
\end{proof}
%=======================================================================================================
\section{Conclusion and Future Work}
In conclusion, this paper has embarked on an extensive investigation into the domain of matrix means interpolation and comparison. A key aspect of this research has been the expansion of the parameter $\vartheta$ from the closed interval $[0,1]$ to encompass the entire positive real line, represented as $\mathbb{R}^+$. This extension has allowed us to explore a broader spectrum of mathematical relationships and properties within this framework.

Furthermore, our exploration has led to the development of various novel results related to Heinz means. We have introduced scalar variants of Heinz inequalities, leveraging Kantorovich's constant, and have extended these inequalities to the operator realm. This expansion not only deepens our understanding of Heinz means but also opens up new avenues for applications in diverse mathematical contexts.

Lastly, we have presented refined Young's type inequalities specifically tailored for traces, determinants, and norms of positive semi-definite matrices. These refined inequalities are expected to find utility in various matrix analysis and linear algebra problems, enhancing our ability to derive meaningful conclusions and insights from the study of positive semi-definite matrices.

As for future work, there are several intriguing directions to consider. Firstly, it may be valuable to explore further extensions of the parameter space beyond $\mathbb{R}^+$ and investigate the implications of such extensions on matrix means and related inequalities. Additionally, the applicability of the developed results in practical fields such as physics, engineering, and data science warrants investigation. Finally, refining and expanding upon the presented inequalities could lead to even more powerful tools for matrix analysis and optimization, offering new insights and solutions to complex problems in mathematics and its applications.

{\bf Conflicts of Interest:} The authors declare that there are no conflicts of interest regarding the
publication of this paper.
%===================================================================================

%===========================================================================
%%%%%%%%%%%%%%%%%%%%%%%%%%%%%%%%%%%%%%%%%%%%%%%%%%%%%%%%%%%%%%%%%%%%%%%%%%%%%%%%%%%%%%%%%%%%

\bibliographystyle{unsrtnat}
\bibliography{references}  %%% Uncomment this line and comment out the ``thebibliography'' section below to use the external .bib file (using bibtex) .




% 	\bibitem{keshet2016prediction}
% 	Keshet, Renato, Alina Maor, and George Kour.
% 	\newblock Prediction-Based, Prioritized Market-Share Insight Extraction.
% 	\newblock In {\em Advanced Data Mining and Applications (ADMA), 2016 12th International
%                       Conference of}, pages 81--94,2016.

% \end{thebibliography}

\end{document}